\documentclass[12pt]{article}
\newtheorem{lemma}{Lemma}[section]
\newtheorem{theorem}[lemma]{Theorem}
\newtheorem{proposition}[lemma]{Proposition}

\newtheorem{definition}[lemma]{Definition}
\newenvironment{proof}{{\bf Proof}}{{\hfill $ \Box $}\vskip 4mm}
\newenvironment{remark}{\addtocounter{lemma}{1}
{\bf Remark \thelemma}}{{\hfill}\vskip 4mm}
\newenvironment{remarks}{\addtocounter{lemma}{1}
{\bf Remarks \thelemma}}{{\hfill}\vskip 4mm}

\newenvironment{examples}{\addtocounter{lemma}{1}
{\bf Examples \thelemma}}{{\hfill}\vskip 4mm}
\newcommand{\nc}{\newcommand}
\nc{\rnc}{\renewcommand}
\nc{\nt}{\newtheorem}

        {~\hfill~\fbox{}\end{trivlist}}

\nc{\thlabel}[1]{\label{theo:#1}}
\nc{\thref}[1]{Theorem~\ref{theo:#1}}
\nc{\selabel}[1]{\label{sect:#1}}
\nc{\seref}[1]{Section~\ref{sect:#1}}
\nc{\lelabel}[1]{\label{lemm:#1}}
\nc{\leref}[1]{Lemma~\ref{lemm:#1}}
\nc{\prlabel}[1]{\label{prop:#1}}
\nc{\prref}[1]{Proposition~\ref{prop:#1}}
\nc{\colabel}[1]{\label{coro:#1}}
\nc{\coref}[1]{Corollary~\ref{coro:#1}}
\nc{\exlabel}[1]{\label{exam:#1}}
\nc{\exref}[1]{Example~\ref{exam:#1}}
\nc{\delabel}[1]{\label{defi:#1}}
\nc{\deref}[1]{Definition~\ref{defi:#1}}
\nc{\eqlabel}[1]{\label{equation:#1}}
\nc{\eqref}[1]{(\ref{equation:#1})}
\nc{\csm}{\mbox{$\triangleright\!\!\!<$}}
\nc{\smc}{\mbox{$>\!\!\!\triangleleft$}}
\nc{\trr}{\triangleright}
\providecommand{\operatorname}[1]{\mathrm{#1}\,}
\nc{\Hom}{\operatorname{Hom}}
\nc{\Mor}{\operatorname{Mor}}
\nc{\Aut}{\operatorname{Aut}}
\nc{\Ann}{\operatorname{Ann}}
\nc{\Ker}{\operatorname{Ker}}
\nc{\Trace}{\operatorname{Trace}}
\nc{\Char}{\operatorname{Char}}
\nc{\Mod}{\operatorname{Mod}}
\nc{\End}{\operatorname{End}}
\nc{\Spec}{\operatorname{Spec}}
\nc{\Span}{\operatorname{Span}}
\nc{\sgn}{\operatorname{sgn}}
\nc{\Id}{\operatorname{Id}}
\nc{\Com}{\operatorname{Com}}

\def\Box{\mbox{$\sqcap\!\!\!\!\sqcup$}}

\nc{\dht}{\mbox{$\rightharpoonup\hspace{-2ex}\rightharpoonup$}}
\nc{\dhtb}{\mbox{$\leftharpoonup\hspace{-2ex}\leftharpoonup$}}
\nc{\nd}{\mbox{$\not|$}} %not divide sign
\providecommand{\text}[1]{\mbox{{\textrm #1}}}

\nc{\nci}{\mbox{$\not\subseteq$}}
\nc{\scontainin}{\mbox{$\mbox{}\subseteq\hspace{-1.5ex}\raisebox{-.5ex}{$_\prime
$}\hspace*{1.5ex}$}}

\def\ot{\otimes}

\def\doublerightleft#1#2{{\lower.2ex\vbox{
\hbox{${\smash{\mathop{\longrightarrow}\limits^{#1}}}$}\vspace*{-4mm}
\hbox{${\smash{\mathop{\longleftarrow}\limits_{#2}}}$}}}}

\def\nint{\hbox{\bbb Z}}

\def\nrea{\hbox{\bbb R}}

\newfont{\bbb}{msbm10 scaled\magstep1}  % Blackboardbold for 12pt article
\newfont{\bbbsub}{msbm10}                % Blackboardbold for subscripts for                                     % 12pt article
\newfont{\msam}{msam10 scaled\magstep1}

\begin{document}
\title{New types of bialgebras arising from the Hopf equation}
\author{G. Militaru
\\University of Bucharest
\\Faculty of Mathematics\\Str. Academiei 14
\\RO-70109 Bucharest 1, Romania
\\e-mail: gmilit@al.math.unibuc.ro}
\date{}
\maketitle
\begin{abstract}
\noindent
Let $M$ be a $k$-vector space and $R\in \End_k(M\ot M)$.
In \cite{M1} we introduced and studied what we called the Hopf equation:
$R^{12}R^{23}=R^{23}R^{13}R^{12}$. By means of a FRT type theorem, we
have proven that the category ${}_H{\cal M}^H$ of $H$-Hopf modules is
deeply involved in solving this equation. In the present paper, we
continue to study the Hopf equation from another perspective: having in
mind the quantum Yang-Baxter equation, in the solution of which the
co-quasitriangular (or braided) bialgebras play an important role (see
\cite{K}), we introduce and study what we call bialgebras with Hopf
functions. The main theorem of this paper shows that, if $M$ is finite
dimensional, any solution $R$ of the Hopf equation has the form
$R=R_{\sigma}$, where $M$ is a right comodule over a bialgebra with a
Hopf function $(B(R),C,\sigma)$ and $R_{\sigma}$ is the special map
$R_{\sigma}(m\ot n)=\sum \sigma(m_{<1>}\ot n_{<1>})m_{<0>}\ot n_{<0>}$.
\end{abstract}

\section{Introduction}
Let $H$ be a Hopf algebra over a field $k$. The strong link between
the category ${}_H{\cal M}^H$ of Hopf modules and the category
${}_H{\cal YD}^H$ of Yetter-Drinfel'd modules recently highlighted
in \cite{CMZ1} (namely, the fact that both are particular cases of
the same general category
${}_A{\cal M}(H)^C$ of Doi-Hopf modules, defined by Doi in \cite{D})
led us in \cite{CMZ2} and \cite{CMZ3} to study the implication of the
category ${}_H{\cal YD}^H$ in the classic, non-quantic part of Hopf
algebra theory.
In \cite{M1} we called this technique "quantisation".
For example, the theorem 4.2 of \cite{CMZ2}, stating that the forgetful
functor ${}_H{\cal YD}^H\rightarrow {}_H{\cal M}$ is Frobenius if and
only if H is finite dimensional and unimodular, can be viewed as the
"quantum version" of the classical theorem saying that any finite
dimensional Hopf algebra is Frobenius.

In \cite{M1} we start to study the reverse problem (it was called
"dequantisation"): that is, we study the category ${}_H{\cal M}^H$ in
connection with problems which so far were specific solely to the
category ${}_H{\cal YD}^H$. The starting
point is simple; it is enough to remember that
the category ${}_H{\cal YD}^H$ is deeply involved in the quantum
Yang-Baxter equation:
$$
R^{12}R^{13}R^{23}=R^{23}R^{13}R^{12}
$$
where $R\in \End_k(M\ot M)$, $M$ being a $k$-vector space.
We evidence the fact that the category ${}_H{\cal M}^H$
can also be studied in connection with a certain non-linear equation.
We call it the Hopf equation, and it is:
$$
R^{12}R^{23}=R^{23}R^{13}R^{12}
$$
The main result of \cite{M1} is a FRT type theorem which shows that in
the finite dimensional case, any solution $R$ of the Hopf equation has
the form $R=R_{(M,\cdot ,\rho )}$, where  $(M,\cdot ,\rho)$
is an object in ${}_{B(R)}{\cal M}^{B(R)}$, for some bialgebra $B(R)$.
In this paper we shall continue to study the Hopf equation from another
perspective. To begin with, we remind that in the quantum Yang-Baxter
equation another important role is played by the co-quasitriangular
bialgebras. These can be viewed as bialgebras $H$ with a k-bilinear map
$\sigma :H\ot H\to k$ satisfying properties which ensure that the
special map
$$
R_{\sigma}:M\ot M\to M\ot M,\quad
R_{\sigma}(m\ot n)=\sum \sigma(m_{<1>}\ot n_{<1>})m_{<0>}\ot n_{<0>}
$$
is a solution for the quantum Yang-Baxter equation.
Starting from here, we introduce new classes of bialgebras which will play
for the Hopf equation the same role as the co-quasitriangular bialgebras do
for the quantum Yang-Baxter equation. We called them $(H,C,\sigma)$
bialgebras with a Hopf function $\sigma :C\ot H\to k$, where $C$ is a
subcoalgebra of $H$.
The reason why the map $\sigma$ is not defined for the entire $H\ot H$,
but only relative to a subcoalgebra $C$ of $H$ is explained in Remarks 2.2
and 2.9. The main result of this paper is theorem 2.8: if $M$ is a
finite dimensional vector space and $R$ is a solution of the Hopf
equation, then for a special subcoalgebra $C$ of $B(R)$ there exists a
unique Hopf function $\sigma: C\ot B(R)\to k$ such that $R=R_{\sigma}$.
We apply the above results by presenting several examples of Hopf functions
on bialgebras. In the last section, as an appendix, we also introduced the
concept corresponding to quasitriangular bialgebras.
\section{Preliminaries}
Throughout this paper, $k$ will be a field.
All vector spaces, algebras, coalgebras and bialgebras considered
are over $k$. $\ot$ and $\Hom$ will mean $\ot_k$ and $\Hom_k$.
For a coalgebra $C$, we will use Sweedler's $\Sigma$-notation, that is,
$\Delta(c)=\sum c_{(1)}\ot c_{(2)},~(I\ot\Delta)\Delta(c)=
\sum c_{(1)}\ot c_{(2)}\ot c_{(3)}$, etc. We will
also use  Sweedler's notation for right $C$-comodules:
$\rho_M(m)=\sum m_{<0>}\otimes m_{<1>}$, for any $m\in M$ if
$(M,\rho_M)$ is a right $ C$-comodule. ${\cal M}^C$ will be the
category of right $C$-comodules and $C$-colinear maps and
${}_A{\cal M}$ will be the category of left $A$-modules and
$A$-linear maps, if $A$ is a $k$-algebra.

From now on, $H$ will be a bialgebra. An element $T\in H^*$ is called
a right integral on $H$ (see \cite{M}) if
$$
Tf=f(1_H)T
$$
for all $f\in H^*$. This is equivalent to
$$
\sum T(h_{(1)})h_{(2)}=T(h)1_H
$$
for all $h\in H$.
Recall that a (left-right) $H$-Hopf module is a
left $H$-module $(M, \cdot)$ which is also a right $H$-comodule $(M,\rho)$
such that
\begin{equation}\label{H1}
\rho(h\cdot m)=\sum h_{(1)}\cdot m_{<0>}\ot h_{(2)}m_{<1>}
\end{equation}
for all $h\in H$, $m\in M$. ${}_H{\cal M}^H$ will be the category of
$H$-Hopf modules and $H$-linear $H$-colinear homomorphisms.
If $(M,\cdot,\rho)\in {}_H{\cal M}^H$ we can define the special map
$$
R_{(M,\cdot ,\rho )}:M\ot M\to M\ot M, \quad
R_{(M,\cdot ,\rho )}(m\ot n)=\sum n_{<1>}\cdot m\ot n_{<0>}.
$$
For a vector space $V$, $\tau :V\otimes V\to V\otimes V$
will denote the switch map, that is, $\tau (v\otimes w)=w\otimes v$
for all $v,w \in V$. If $R:V\ot V\to V\ot V$ is a linear map
we denote by $R^{12}$, $R^{13}$, $R^{23}$ the maps of $\End_k(V\ot V\ot V)$
given by
$$
R^{12}=R\ot I, \quad R^{23}=I\ot R,\quad
R^{13}=(I\ot \tau)(R\ot I)(I\ot \tau).
$$
Using the notation $R(u\ot v)=\sum u_1\ot v_1$, then
$$
R^{12}(u\ot v\ot w)=\sum u_1\ot v_1\ot w_0
$$
where the subscript (0) means that $w$ is not affected by the application
of $R^{12}$.\\
Let $H$ be a bialgebra and $(M, \cdot)$ a left $H$-module which is also a
right $H$-comodule $(M,\rho)$. We recall that $M(,\cdot,\rho)$ is a
Yetter-Drinfel'd module if the following compatibility relation holds:
$$
\sum h_{(1)}\cdot m_{<0>}\ot h_{(2)}m_{<1>}
=\sum (h_{(2)}\cdot m)_{<0>}\ot (h_{(2)}\cdot m)_{<1>}h_{(1)}
$$
for all $h\in H$, $m\in M$. ${}_H{\cal YD}^H$ will be the category of
Yetter-Drinfel'd modules and $H$-linear $H$-colinear homomorphism.
For a further study of the Yetter-Drinfel'd category we refer to
\cite{LR}, \cite{R1}, \cite{RT}, \cite{Y}, or to the more recent
\cite{CMZ1}, \cite{CMZ2}, \cite{CMZ3}, \cite{FMS}.

Let $H$ be a bialgebra over $k$ and $\sigma :H\ot H\to k$ be a
$k$-bilinear map. Recall that the pair $(H,\sigma)$ is a
{\sl co-quasitriangular} (or {\sl braided}) bialgebra if

$(B1)\quad \sum \sigma(x_{(1)}\ot y_{(1)})y_{(2)}x_{(2)}=
\sum \sigma(x_{(2)}\ot y_{(2)})x_{(1)}y_{(1)}$

$(B2)\quad \sigma (x\ot 1)=\varepsilon(x)$

$(B3) \quad \sigma(x\ot yz)=\sum \sigma (x_{(1)}\ot y)\sigma (x_{(2)}\ot z)$

$(B4)\quad \sigma (1\ot x)=\varepsilon(x)$

$(B5) \quad \sigma(xy\ot z)=\sum \sigma (y\ot z_{(1)})\sigma (x\ot z_{(2)})$

for all $x$, $y$, $z\in H$. If $(H,\sigma)$ is a co-quasitriangular
bialgebra and $(M,\rho)$ is a right $H$-comodule, then the special
map
$$
R_{\sigma}:M\ot M\to M\ot M, \quad
R_{\sigma}(m\ot n)=\sum \sigma(m_{<1>}\ot n_{<1>})m_{<0>}\ot n_{<0>}
$$
is a solution for the quantum Yang-Baxter equation
$$
R^{12}R^{13}R^{23}=R^{23}R^{13}R^{12}.
$$
Conversely, if $M$ is a finite dimensional vector space and $R$ is a
solution of the quantum Yang-Baxter equation, then there exists a
bialgebra $A(R)$ and a unique $k$-bilinear map $\sigma :A(R)\ot A(R)\to k$
such that $(A(R), \sigma)$ is co-quasitriangular,
$M\in {\cal M}^{A(R)}$ and $R=R_{\sigma}$.

A quasitriangular bialgebra is a pair $(H,R)$, where $H$ is a
bialgebra and $R\in H\otimes H$ such that the following conditions are
fulfilled:

$(QT1) \quad \sum \Delta(R^1)\otimes R^2=R^{13}R^{23}$

$(QT2) \quad \sum \varepsilon(R^1)R^2=1$

$(QT3) \quad \sum R^1\otimes \Delta(R^2)=R^{13}R^{12}$

$(QT4) \quad \sum R^1\varepsilon(R^2)=1$

$(QT5) \quad \Delta^{\rm cop}(h)R=R\Delta(h)$, for all $h\in H$.

Recall from \cite{M1} the following:
\begin{definition}
Let $M$ be a vector space and $R\in \End_k(M\ot M)$.
We say that $R$ is a solution for the Hopf equation if
\begin{equation}\label{Heq}
R^{23}R^{13}R^{12}=R^{12}R^{23}
\end{equation}
\end{definition}
The following lemma will be important
\begin{lemma}\label{hei}
Let $M$ be a finite dimensional vector space and $\{m_1,\cdots,m_n \}$
a basis of $M$. Let $R$, $S\in \End_k(M\ot M)$ given by
$$
R(m_v\ot m_u)=\sum_{i,j}x_{uv}^{ji}m_i\ot m_j, \quad
S(m_v\ot m_u)=\sum_{i,j}y_{uv}^{ji}m_i\ot m_j,
$$
for all $u$, $v=1,\cdots ,n$, where $(x_{uv}^{ji})_{i,j,u,v}$,
$(y_{uv}^{ji})_{i,j,u,v}$ are two families of scalars of $k$. Then
$$
R^{23}S^{13}S^{12}=S^{12}R^{23}
$$
if and only if
$$
\sum_{u,v,\beta}x_{uv}^{ji}y_{k\beta}^{up}y_{lq}^{v\beta}=
\sum_{\alpha}x_{kl}^{j\alpha}y_{\alpha q}^{ip}
$$
for all $i$, $j$, $k$, $l$, $p$, $q=1,\cdots, n$. In particular, $R$ is a
solution for the Hopf equation if and only if
\begin{equation}\label{tam}
\sum_{u,v,\beta}x_{uv}^{ji}x_{k\beta}^{up}x_{lq}^{v\beta}=
\sum_{\alpha}x_{kl}^{j\alpha}x_{\alpha q}^{ip}
\end{equation}
for all $i$, $j$, $k$, $l$, $p$, $q=1,\cdots, n$.
\end{lemma}

\begin{proof}
For $k$, $l$, $q=1,\cdots, n$ we have:
\begin{eqnarray*}
R^{23}S^{13}S^{12}(m_q\ot m_l\ot m_k )&=&
R^{23}S^{13}\Bigl(\sum_{\beta, v}
y_{lq}^{v\beta}m_{\beta}\ot m_v\ot m_{k} \Bigl)\\
&=&R^{23}\Bigl(\sum_{\beta, v,u,p}y_{k\beta}^{up}
y_{lq}^{v\beta}m_{p}\ot m_v\ot m_{u} \Bigl)\\
&=& \sum_{i,j,p} \Bigl ( \sum_{u,v,\beta}x_{uv}^{ji}
y_{k\beta}^{up}y_{lq}^{v\beta} \Bigl ) m_{p}\ot m_i\ot m_j
\end{eqnarray*}
and
\begin{eqnarray*}
S^{12}R^{23}(m_q\ot m_l\ot m_k )&=&
S^{12}\Bigl(\sum_{j,\alpha}x_{kl}^{j\alpha}m_q\ot m_{\alpha}\ot m_j\Bigl)\\
&=&\sum_{j,\alpha, p,i}x_{kl}^{j\alpha}y_{\alpha q}^{ip}
m_p\ot m_{i}\ot m_j\\
&=&\sum_{i,j,p}\Bigl(\sum_{\alpha}x_{kl}^{j\alpha}y_{\alpha q}^{ip} \Bigl)
m_p\ot m_{i}\ot m_j
\end{eqnarray*}
Hence, the conclusion follows.
\end{proof}

Recall now the main results of \cite{M1}.

\begin{theorem}
Let $M$ be a finite dimensional vector space and $R\in \End_k(M\ot M)$
be a solution of the Hopf equation. Then
\begin{enumerate}
\item There exists a bialgebra $B(R)$ such that $M$ has a structure
of $B(R)$-Hopf module $(M,\cdot, \rho)$ and $R=R_{(M,\cdot, \rho)}$.
\item The bialgebra $B(R)$ is a universal object with this property:
if $H$ is a bialgebra such that
$(M,\cdot^{\prime}, \rho^{\prime})\in {}_H{\cal M}^H$ and
$R=R_{(M,\cdot^{\prime}, \rho^{\prime})}$ then there exists a unique
bialgebra map $f:B(R)\to H$ such that $\rho^{\prime}=(I\ot f)\rho$.
Furthermore, $a\cdot m=f(a)\cdot^{\prime}m$, for all $a\in B(R)$,
$m\in M$.
\end{enumerate}
\end{theorem}

\section{Hopf function on a bialgebra}
First of all we will introduce the following key definition of this
paper:
\begin{definition}
Let $H$ be a bialgebra and $C$ be a subcoalgebra of $H$. A $k$-biliniar
map $\sigma :C\ot H\to k$ is called a Hopf function if:

$(H1)\quad \sum \sigma(c_{(1)}\ot h_{(1)})h_{(2)}c_{(2)}=
\sum \sigma(c_{(2)}\ot h)c_{(1)}$

$(H2)\quad \sigma (c\ot 1)=\varepsilon(c)$

$(H3) \quad \sigma(c\ot hk)=\sum \sigma (c_{(1)}\ot h)\sigma (c_{(2)}\ot k)$

for all $c\in C$, $h$, $k\in H$. In this case we shall say that
$(H,C,\sigma)$ is a bialgebra with a Hopf function.
\end{definition}

\begin{remarks}
1. The first question which arises is why we have not defined the map
$\sigma$ on the entire $H\ot H$ and we have instead presented a
definition relative to a subcoalgebra $C$ of $H$. This choice was driven
by the compatibility condition (H1): in the case $\sigma :H\ot H\to k$
and setting $c=1_H$, this would result in the integral type relation
\begin{equation}\label{unua}
\sigma (1_H\ot h)1_H=\sum \sigma(1_H\ot h_{(1)})h_{(2)}
\end{equation}
for all $h\in H$. Hence, from (H3), the map
$T_{\sigma}:H\to k$, $T_{\sigma}(h):=\sigma(1_H\ot h)$, is an algebra
map, and from (\ref{unua}), a right integral on $H$. Using lemma 2.2 of
\cite{MS}, we are led to the trivial $H=k$. This is why defining
$\sigma$ relative to a subcoalgebra of $H$ becomes mandatory.

2. The conditions of compatibility (H2) and (H3) are exactly
(B2) and (B3), respecting the definition relativ to $C$. The left
hand side of (H1) is the same with the left hand side of (B1),
while the right hand side has suffered, as we expect, considerable changes.

3. Let $(H,C,\sigma)$ be a bialgebra with a Hopf function. If $\sigma$ is
right invertible in the convolution algebra $\Hom_k(C\ot H, k)$, then
(H2) follows from (H3). Indeed, for $c\in C$ we have:
\begin{eqnarray*}
\sigma (c\ot 1)&=&\sum \sigma (c_{(1)}\ot 1)\varepsilon (c_{(2)})\\
&=&\sum \sigma (c_{(1)}\ot 1)\sigma (c_{(2)}\ot 1)\sigma^{-1}(c_{(3)}\ot 1)\\
&=&\sum \sigma (c_{(1)}\ot 1)\sigma^{-1}(c_{(2)}\ot 1)=\varepsilon(c)
\end{eqnarray*}
If $H$ has an antipode $S$, then $\sigma$ is invertible and
$\sigma^{-1}(c\ot h)=\sigma (c\ot S(h))$, for all $c\in C$, $h\in H$.
\end{remarks}

We shall point out the link existing between the (H1) compatibility
condition and the concept of right integral on $H$.

Let $H$ be a bialgebra and $C$ a subcoalgebra of $H$. If $T\in H^*$ is a
right integral on $H$ then the map
$$
\sigma_{T}:C\ot H\to k,\qquad
\sigma_{T}(c\ot h):=\varepsilon(c)T(h), \quad \forall c\in C, h\in H
$$
satisfies (H1).

Conversely, if $1_H\in C$ and $\sigma:C\ot H\to k$ satisfies (H1)
then the map
$$
T_{\sigma}:H\to k,\qquad T_{\sigma}(h):=\sigma(1_H\ot h),\quad
\forall h\in H
$$
is a right integral on $H$. In addition, we suppose that $H$ has an
antipode and (H2) also holds. Then, $T_{\sigma}(1_H)=1_k$; so, using the
classical dual Maschke theorem for Hopf algebras (see \cite{M}) we obtain
that $H$ is cosemisimple. As $T_{\sigma_{T}}=T$, we obtain that the map
$$
\{\sigma:C\ot H\to k \mid \sigma \;\mbox{satisfies}\; (H1) \}\to
\int_{H^{*}}^{r}, \quad \sigma\to T_{\sigma}
$$
is surjective, and the map
$$
\int_{H^{*}}^{r}\to \{\sigma:C\ot H\to k
\mid \sigma\; \mbox{satisfies}\; (H1) \}, \quad T\to \sigma_{T}
$$
is injective. We record these observation in the following:

\begin{proposition}
Let $H$ be a bialgebra and $C$ a subcoalgebra of $H$. Then:
\begin{enumerate}
\item if $T\in H^*$ is a right integral on $H$, then the map
$\sigma_{T}:C\ot H\to k$ satisfies (H1).
\item if $1_H\in C$ and $\sigma:C\ot H\to k$ satisfies (H1),
then the map $T_{\sigma}:H\to k$ is a right integral on $H$.
Furthermore, if (H2) holds and $H$ has an antipode, then $H$ is
cosemisimple.
\end{enumerate}
\end{proposition}

In the next proposition we shall prove that, if a $k$-bilinear map
$\sigma :C\ot H\to k$ satisfies (H3) and (H1) holds for a basis of $C$
and a sistem of generators of $H$ as an algebra, then (H1) holds for any
$c\in C$ and $h\in H$.

\begin{proposition}\label{gen}
Let $H$ be a bialgebra, $C$ be a subcoalgebra of $H$ and
$\sigma :C\ot H\to k$ a $k$-bilinear map which satisfies (H3).
Suppose that (H1) holds for a basis of $C$ and a sistem of generators
of $H$ as an algebra. Then (H1) holds for any $c\in C$ and $h\in H$.
\end{proposition}

\begin{proof}
Let $c\in C$ be an element of the given basis and $x$, $y\in H$ two
elements between the generators of $H$. It is enough to prove that
(H1) holds for $(c,xy)$. We have:
\begin{eqnarray*}
\sum \sigma (c_{(2)}\ot xy)c_{(1)}&=&
\sum \sigma (c_{(2)(1)}\ot x) \sigma (c_{(2)(2)}\ot y)c_{(1)}\\
&=&\sum \sigma (c_{(2)}\ot y) \sigma (c_{(1)(2)}\ot x)c_{(1)(1)}\\
\text{(\mbox{(H1) holds for} $x$ )}
&=&\sum\sigma (c_{(2)}\ot y)\sigma (c_{(1)(1)}\ot x_{(1)})x_{(2)}c_{(1)(2)}\\
&=&\sum\sigma (c_{(1)}\ot x_{(1)})x_{(2)}\sigma (c_{(2)(2)}\ot y)c_{(2)(1)}\\
\text{(\mbox{(H1) holds for} $y$ )}
&=&\sum \sigma (c_{(1)}\ot x_{(1)})x_{(2)}\sigma (c_{(2)(1)}\ot y_{(1)})
y_{(2)}c_{(2)(2)}\\
&=&\sum \sigma (c_{(1)}\ot x_{(1)})\sigma (c_{(2)}\ot y_{(1)})x_{(2)}
y_{(2)}c_{(3)}\\
\text{(\mbox{using (H3)} )}
&=&\sum \sigma (c_{(1)}\ot x_{(1)}y_{(1)})x_{(2)}y_{(2)}c_{(2)}\\
&=&\sum \sigma (c_{(1)}\ot (xy)_{(1)})(xy)_{(2)}c_{(2)}
\end{eqnarray*}
and we are done.
\end{proof}

To give examples of Hopf functions on the classical examples of bialgebras
sems to be a difficult problem (see examples 1) and 2) below). More examples
will be given after the main results of this paper.

\begin{examples} \label{25}
1. Let $G$ be a nontrivial group and $H=k[G]$ be the groupal Hopf
algebra. Let $C$ be an arbitrary subcoalgebra of $H$. Then there exists
no Hopf function $\sigma :C\ot k[G]\to k$.

Indeed, any subcoalgebra of $k[G]$ has the form $k[F]$, where $F$ is a
subset of $G$. Suppose that there exists $\sigma :k[F]\ot k[G]\to k$
a Hopf function. From (H2) we get that $\sigma (f\ot 1)=1$ for all
$f\in F$. Now let $g\in G$, $g\neq 1$ and $f\in F$. From (H1) we obtain
$\sigma(f\ot g)fg=\sigma(f\ot g)f$ i.e. $\sigma(f\ot g)=0$ for all
$f\in F$. But then, using (H3), we get
$$
1=\sigma(f\ot 1)=\sigma(f\ot gg^{-1})=\sigma(f\ot g)\sigma(f\ot g^{-1})=0,
$$
contradiction.

2. Let $H=k[X,X^{-1}]$, which is a bialgebra with $\Delta(X)=X\ot X$,
$\varepsilon(X)=1$. Let $C$ be an arbitrary subcoalgebra of $H$.
Then there exists no Hopf function $\sigma :C\ot H\to k$.

Let us suppose that $C$ is a subcoalgebra of $H$ and $\sigma :C\ot H\to k$
is a Hopf function. Then there exists $t\in \nint^{*}$ such that
$X^t\in C$. Then, by (H2), $\sigma(X^t\ot 1)=1$ and from (H1) we obtain
that $\sigma(X^t\ot X)X^{t+1}=\sigma(X^t\ot X)X^{t}$, i.e.
$\sigma(X^t\ot X)=0$. But then, using (H3) we get
$$
1=\sigma(X^t\ot 1)=\sigma(X^t\ot XX^{-1})=
\sigma(X^t\ot X)\sigma(X^t\ot X^{-1})=0,
$$
contradiction.

3. Let $H=k_{q}<x,y\mid xy=qyx>$ be the quantum plane:
$$
\Delta(x)=x\ot x, \quad \Delta(y)=y\ot 1+x\ot y, \quad
\varepsilon(x)=1, \quad \varepsilon(y)=0.
$$
Let $C:=kx$ be the one dimensional subcoalgebra of $H$ with $\{x\}$ a
$k$-basis. Let $a\in k$ and $\sigma_{a}:C\ot H\to k$ given by
$$
\sigma_{a}(x\ot 1)=1,\quad \sigma_{a}(x\ot x)=0,\quad
\sigma_{a}(x\ot y)=a,
$$
and extend $\sigma_{a}$ to the entire $C\ot H$ with (H3). Then, $\sigma_{a}$
is a Hopf function.

Using proposition (\ref{gen}), it is enough to check that (H1) holds for
$h\in \{x,y\}$. For $h=x$, (H1) is
$$
\sigma_{a}(x\ot x)x^2=\sigma_{a}(x\ot x)x
$$
which holds if and only if $\sigma_{a}(x\ot x)=0$. For $h=y$, (H1) has
the form
$$
\sigma_{a}(x\ot y)x+\sigma_{a}(x\ot x)yx=\sigma_{a}(x\ot y)x,
$$
which is true, as $\sigma_{a}(x\ot x)=0$. In fact, we also prove the
converse: if $\sigma:C\ot H\to k$ is a Hopf function, then there
exists $a\in k$ such that $\sigma=\sigma_{a}$.

4. Let ${\cal T}(k)$ be the three dimensional noncocommutative bialgebra
constructed in \cite{M1}, i.e.

$\bullet$ As a vector space, ${\cal T}(k)$ is three dimensional
with $\{1,x,z \}$ a $k$-basis.

$\bullet$ The multiplication rule is given by:
$$
x^2=x, \quad xz=zx=z^2=0.
$$
$\bullet$ The comultiplication $\Delta$ and the counity $\varepsilon$
are given by
$$
\Delta(x)=x\ot x, \quad \Delta(z)=x\ot z+z\ot 1,\quad
\varepsilon(x)=1, \quad \varepsilon(z)=0.
$$
Let $C:=kx$ be the one dimensional subcoalgebra of ${\cal T}(k)$ with
$\{x\}$ a $k$-basis. Then it is easy to see that the $k$-bilinear map
$$
\sigma:C\ot {\cal T}(k)\to k, \quad
\sigma(x\ot 1)=\sigma(x\ot x)=1,\quad \sigma(x\ot z)=0,
$$
defines a Hopf function.

5. Let $M$ be a monoid and $N:=\{n\in M\mid xn=n, \forall x\in M \}$.
Let $H=k[M]$ and $C:=k[F]$, where $F$ is a subset of $N$. Let
$\sigma :k[F]\ot k[M]\to k$ such that
\begin{enumerate}
\item $\sigma(f\ot 1)=1$
\item $\sigma(f,\bullet):M\to (k,\cdot)$ is a morphism of monoids
\end{enumerate}
for all $f\in F$. Then $\sigma$ is a Hopf function.

We shall give such an example. Let $a \in k$ and
${\cal F}_{a}(k)=\{u:k\to k\mid u(a)=a \}$ be the monoid (with the
usual composition of functions), of all function $u$ on $k$ such that
$a$ is a fixed point of $u$. Let $F=\{f\}$, where $f:k\to k$ is the constant
function $f(x)=a$ for all $x\in k$. Then
$$
\sigma :k[F]\ot k[{\cal F}_{a}(k)]\to k, \quad
\sigma(f\ot u)=1
$$
for all $u\in {\cal F}_{a}(k)$, is a Hopf function.
\end{examples}

Let $H$ be a bialgebra, $C$ a subcoalgebra of $H$ and $\sigma :C\ot H\to k$
a $k$-bilinear map. We denote by $\sigma_{12}$, $\sigma_{13}$,
$\sigma_{23}$ the maps of $\Hom_k(C\ot C\ot H, k)$ given by:
$$
\sigma_{12}(c\ot d\ot x):=\varepsilon(x)\sigma(c\ot d),\quad
\sigma_{13}(c\ot d\ot x):=\varepsilon(d)\sigma(c\ot x),\quad
\sigma_{23}(c\ot d\ot x):=\varepsilon(c)\sigma(d\ot x)
$$
for all $c$, $d\in C$, $x\in H$.

\begin{proposition}
Let $(H,C,\sigma)$ be a bialgebra with a Hopf function $\sigma :C\ot H\to k$.
Then:
\begin{enumerate}
\item in the convolution algebra $\Hom_k(C\ot C\ot H, k)$ the following
identity holds:
\begin{equation}\label{dec}
\sigma_{23}*\sigma_{13}*\sigma_{12}=\sigma_{12}*\sigma_{23}
\end{equation}
\item if $(M,\rho)$ is a right $C$-comodule, then the special map
$$
R_{\sigma}:M\ot M\to M\ot M, \quad
R_{\sigma}(m\ot n)=\sum \sigma(m_{<1>}\ot n_{<1>})m_{<0>}\ot n_{<0>}
$$
is a solution for the Hopf equation.
\end{enumerate}
\end{proposition}

\begin{proof}
1. Let $c$, $d\in C$ and $x\in H$. We have:
\begin{eqnarray*}
(\sigma_{12}*\sigma_{23})(c\ot d\ot x)&=&
\sum \varepsilon(x_{(1)})\sigma(c_{(1)}\ot d_{(1)})
\varepsilon(c_{(2)})\sigma(d_{(2)}\ot x_{(2)})\\
&=&\sum \sigma(c\ot d_{(1)})\sigma(d_{(2)}\ot x)
\end{eqnarray*}
and

$(\sigma_{23}*\sigma_{13}*\sigma_{12})(c\ot d\ot x)=$
\begin{eqnarray*}
&=&\sum \varepsilon(c_{(1)})\sigma(d_{(1)}\ot x_{(1)})
\varepsilon(d_{(2)})\sigma(c_{(2)}\ot x_{(2)})
\varepsilon(x_{(3)})\sigma(c_{(3)}\ot d_{(3)})\\
&=&\sum\underline{\sigma (c_{(1)}\ot x_{(2)})\sigma (c_{(2)}\ot d_{(2)})}
\sigma (d_{(1)}\ot x_{(1)})\\
\text{(\mbox{using (H3)})}
&=&\sum \sigma(c\ot x_{(2)}d_{(2)})\sigma (d_{(1)}\ot x_{(1)})\\
&=&\sum \sigma
\Bigl(c\ot\underline{\sigma(d_{(1)}\ot x_{(1)})x_{(2)}d_{(2)}} \Bigl)\\
\text{(\mbox{using (H1)})}
&=&\sum \sigma
\Bigl(c\ot \sigma(d_{(2)}\ot x )d_{(1)} \Bigl)\\
&=&\sum \sigma(c\ot d_{(1)})\sigma(d_{(2)}\ot x)
\end{eqnarray*}
i.e. the formula (\ref{dec}) holds.

2. Let $R=R_{\sigma}$ and $u$, $v$, $w\in M$. Then, the fact that $R$ is a
solution of the Hopf equation will follow from equation (\ref{dec})
and from the formulas:
$$
R^{12}R^{23}(u\ot v\ot w)=\sum
\Bigl(\sigma_{12}*\sigma_{23}\Bigl)(u_{<1>}\ot v_{<1>}\ot w_{<1>})
u_{<0>}\ot v_{<0>}\ot w_{<0>}
$$
and
$$
R^{23}R^{13}R^{12}(u\ot v\ot w)=\sum
\Bigl(\sigma_{23}*\sigma_{13}*\sigma_{12}\Bigl)
(u_{<1>}\ot v_{<1>}\ot w_{<1>})u_{<0>}\ot v_{<0>}\ot w_{<0>}
$$
Indeed, we have:
\begin{eqnarray*}
R^{12}R^{23}(u\ot v\ot w)&=&R^{12}
\Bigl(\sum\sigma (v_{<1>}\ot w_{<1>})u\ot v_{<0>}\ot w_{<0>}\Bigl)\\
&=&\sum\sigma (u_{<1>}\ot v_{<0><1>})\sigma (v_{<1>}\ot w_{<1>})
u_{<0>}\ot v_{<0><0>}\ot w_{<0>}\\
&=&\sum\sigma (u_{<1>}\ot v_{<1>(1)})\sigma (v_{<1>(2)}\ot w_{<1>})
u_{<0>}\ot v_{<0>}\ot w_{<0>}\\
&=&\sum\Bigl(\sigma_{12}*\sigma_{23}\Bigl)(u_{<1>}\ot v_{<1>}\ot w_{<1>})
u_{<0>}\ot v_{<0>}\ot w_{<0>}
\end{eqnarray*}
On the other hand

$R^{23}R^{13}R^{12}(u\ot v\ot w)=$
\begin{eqnarray*}
&=&R^{23}R^{13}
\Bigl(\sum\sigma (u_{<1>}\ot v_{<1>}) u_{<0>}\ot v_{<0>}\ot w \Bigl)\\
&=&R^{23}
\Bigl(\sum\sigma (u_{<0><1>}\ot w_{<1>})\sigma (u_{<1>}\ot v_{<1>})
u_{<0><0>}\ot v_{<0>}\ot w_{<0>} \Bigl)\\
&=&R^{23}
\Bigl(\sum\sigma (u_{<1>}\ot w_{<1>})\sigma (u_{<2>}\ot v_{<1>})
u_{<0>}\ot v_{<0>}\ot w_{<0>} \Bigl)\\
&=&\sum \sigma (v_{<0><1>}\ot w_{<0><1>})
\sigma (u_{<1>}\ot w_{<1>})\sigma (u_{<2>}\ot v_{<1>})
u_{<0>}\ot v_{<0><0>}\ot w_{<0><0>}\\
&=&\sum \sigma (v_{<1>}\ot w_{<1>})
\sigma (u_{<1>}\ot w_{<2>})\sigma (u_{<2>}\ot v_{<2>})
u_{<0>}\ot v_{<0>}\ot w_{<0>}\\
&=&\sum\sigma (u_{<1>(1)}\ot w_{<1>(2)})\sigma (u_{<1>(2)}\ot v_{<1>(2)})
\sigma (v_{<1>(1)}\ot w_{<1>(1)})
u_{<0>}\ot v_{<0>}\ot w_{<0>}\\
&=&\sum\Bigl(\sigma_{23}*\sigma_{13}*\sigma_{12}\Bigl)
(u_{<1>}\ot v_{<1>}\ot w_{<1>})u_{<0>}\ot v_{<0>}\ot w_{<0>}
\end{eqnarray*}
and the proof is complete now.
\end{proof}

\begin{remark}
There is another proof of the second statement of the above proposition.
Let $(H,C,\sigma)$ be a bialgebra with a Hopf function
$\sigma :C\ot H\to k$. Then, we can construct a functor
$$
F_{\sigma}:{\cal M}^C\to {}_H{\cal M}^H
$$
given as follows: if $(M,\rho)$ is a right $C$-comodule then
$F_{\sigma}(M):=M$ has a structure of right $H$-comodule via
$$
M\stackrel{\rho}{\longrightarrow}
M\ot C\stackrel{I\ot i}{\longrightarrow}M\ot H
$$
where $i:C\to H$ is the inclusion; $M$ has also a left $H$-module
structure given by
$$
h\cdot m:=\sum \sigma(m_{<1>}\ot h)m_{<0>}
$$
for all $h\in H$ and $m\in M$. Furthermore,
$(M,\cdot,\rho)\in {}_H{\cal M}^H$ and $R_{\sigma}=R_{(M,\cdot,\rho)}$.
Now the fact that $R_{\sigma}$ is a solution for the Hopf equation
follows from proposition (2.6) of \cite{M1}.

Indeed, (H2) and (H3) give us that $(M,\cdot)$ is a left $H$-module.
We shall prove that $(M,\cdot,\rho)$ is an $H$-Hopf module using (H1).
For $h\in H$ and $m\in M$ we have:
\begin{eqnarray*}
\sum h_{(1)}\cdot m_{<0>}\ot h_{(2)}m_{<1>}&=&
\sum \sigma (m_{<0><1>}\ot h_{(1)} ) m_{<0><0>}\ot h_{(2)}m_{<1>}\\
&=&\sum m_{<0>} \ot
\underline{\sigma (m_{<1>(1)} \ot h_{(1)} ) h_{(2)}m_{<1>(2)} }\\
\text{(\mbox{using (H1)} )}
&=&\sum m_{<0>}\ot \sigma(m_{<1>(2)}\ot h)m_{<1>(1)}\\
&=& \sum \sigma (m_{<1>}\ot h)m_{<0><0>}\ot m_{<0><1>}\\
&=&\rho \Bigl(\sum \sigma(m_{<1>}\ot h)m_{<0>}\Bigl)\\
&=&\rho(h\cdot m)
\end{eqnarray*}
On the other hand
\begin{eqnarray*}
R_{(M,\cdot, \rho)}(m\ot n)&=&\sum n_{<1>}\cdot m\ot n_{<0>}\\
&=&\sum \sigma (m_{<1>}\ot n_{<1>})m_{<0>}\ot n_{<0>}\\
&=&R_{\sigma}(m\ot n)
\end{eqnarray*}
for all $m$, $n\in M$. We conclude that
$(M,\cdot,\rho)\in {}_H{\cal M}^H$ and $R_{\sigma}=R_{(M,\cdot,\rho)}$.
\end{remark}

We recall from \cite{M1} that if $M$ is a finite dimensional vector space
and $R\in \End_k(M\ot M)$ is a solution for the Hopf equation, then
there exists a bialgebra $B(R)$ such that $M$ has a structure
of $B(R)$-Hopf module $(M,\cdot, \rho)$ and $R=R_{(M,\cdot, \rho)}$.
We recall the construction of $B(R)$:

Let $\{m_1,\cdots ,m_n \}$
be a basis for $M$ and $(x_{uv}^{ji})_{i,j,u,v}$ a family
of scalars of $k$ such that
\begin{equation}
R(m_v\ot m_u)=\sum_{i,j}x_{uv}^{ji}m_i\ot m_j
\end{equation}
for all $u$, $v=1,\cdots ,n$.

Let $(C, \Delta, \varepsilon)={\cal M}^n(k)$, be the comatrix coalgebra
of order $n$, i.e. $C$ is the coalgebra with the basis
$\{c_{ij}\mid i,j=1,\cdots,n \}$ such that
\begin{equation}\label{com1}
\Delta(c_{jk})=\sum_{u=1}^{n}c_{ju}\ot c_{uk},
\quad \varepsilon (c_{jk})=\delta_{jk}
\end{equation}
for all $j,k=1,\cdots, n$. Then $B(R)$ is the free algebra generated
by $(c_{ij})$ with the relations:
$$
\chi(i,j,k,l)=0
$$
where
\begin{equation}\label{obst}
\chi(i,j,k,l):=\sum_{u,v}x_{uv}^{ji}c_{uk}c_{vl} -
\sum_{\alpha}x_{kl}^{j\alpha}c_{i\alpha}
\end{equation}
for all $i$, $j$, $k$, $l=1,\cdots, n$.
$M$ has a right $B(R)$-comodule structure which extends the right
$C$-comodule structure given by
\begin{equation}\label{ro}
\rho(m_l)=\sum_{v=1}^{n}m_v\ot c_{vl}
\end{equation}
for all $l=1,\cdots, n$.

We shall prove now the main result of the paper.

\begin{theorem}
Let $M$ be a finite dimensional vector space and $R\in \End_k(M\ot M)$
be a solution for the Hopf equation. Let $C$ be the subcoalgebra of $B(R)$
with $(c_{ij})$ a $k$-sistem of generators of $C$. Then:
\begin{enumerate}
\item There exists a unique Hopf function
$\sigma:C\ot B(R)\to k$ such that  $R=R_{\sigma}$.
\item If $R$ is bijective and $R^{12}R^{13}=R^{13}R^{12}$, then $\sigma$
is invertible in the convolution algebra $\Hom_k(C\ot B(R),k)$.
\end{enumerate}
\end{theorem}

\begin{proof}
1. First we prove the uniqueness. Let $\sigma:C\ot B(R)\to k$ be a Hopf
function such that  $R=R_{\sigma}$. Let $u$, $v=1,\cdots,n$. Then
\begin{eqnarray*}
R_{\sigma}(m_v\ot m_u)&=&\sum \sigma \Bigl((m_v)_{<1>}\ot (m_u)_{<1>}\Bigl)
(m_v)_{<0>}\ot (m_u)_{<0>}\\
&=&\sum_{i,j} \sigma(c_{iv}\ot c_{ju})m_i\ot m_j
\end{eqnarray*}
Hence $R_{\sigma}(m_v\ot m_u)=R(m_v\ot m_u)$ gives us
\begin{equation}\label{uns}
\sigma(c_{iv}\ot c_{ju})=x_{uv}^{ji}
\end{equation}
for all $i$, $j$, $u$, $v=1,\cdots, n$. As $B(R)$ is generated as an
algebra by $(c_{ij})$, the relations (\ref{uns}) with (H2) and (H3)
ensure the uniqueness of $\sigma$.

Now we shall prove the existence of $\sigma$. First we define
$\sigma_0:C\ot C\to k$ by the formulas (\ref{uns}). Then we extend
$\sigma_{0}$ to
a map $\sigma_1:C\ot T(C)\to k$ such that (H2) and (H3) hold. In order to
prove that $\sigma_1$ factorizes to a map $\sigma:C\ot B(R)\to k$, we have
to show that $\sigma_1(C\ot I)=0$, where $I$ is the two-sided ideal
of $T(C)$ generated by all $\chi(i,j,k,l)$. It is enought to prove that
$$
\sigma_1(c_{pq}\ot \chi(i,j,k,l))=0
$$
for all $i$, $j$, $k$, $l$, $p$, $q=1,\cdots, n$. We have:
\begin{eqnarray*}
\sigma_1(c_{pq}\ot \chi(i,j,k,l))&=&
\sum_{u,v}x_{uv}^{ji}\sigma_1(c_{pq}\ot c_{uk}c_{vl})-
\sum_{\alpha}x_{kl}^{j\alpha}\sigma_1(c_{pq}\ot c_{i\alpha})\\
&=&\sum_{u,v,\beta}x_{uv}^{ji}\sigma_1(c_{p\beta}\ot c_{uk})
\sigma_1(c_{\beta q}\ot c_{vl})-
\sum_{\alpha}x_{kl}^{j\alpha}x_{\alpha q}^{ip}\\
&=&\sum_{u,v,\beta}x_{uv}^{ji}x_{k\beta}^{up}x_{lq}^{v\beta}-
\sum_{\alpha}x_{kl}^{j\alpha}x_{\alpha q}^{ip}\\
\text{(from (\ref{tam}))}
&=&0
\end{eqnarray*}
We conclude that we have constructed $\sigma:C\ot B(R)\to k$ such that
(H2) and (H3) hold and $R=R_{\sigma}$. It remains to prove that (H1)
also holds. Using proposition (\ref{gen}), it is enought to check (H1) for
$c=c_{il}$ and $h=c_{jk}$, for all $i$, $j$, $k$, $l=1,\cdots, n$.
We have:
\begin{eqnarray*}
\sum \sigma(c_{(1)}\ot h_{(1)})h_{(2)}c_{(2)}&=&
\sum_{u,v}\sigma(c_{iv}\ot c_{ju})c_{uk}c_{vl}\\
&=&\sum_{u,v}x_{uv}^{ji}c_{uk}c_{vl}
\end{eqnarray*}
and
\begin{eqnarray*}
\sum \sigma(c_{(2)}\ot h)hc_{(1)}&=&
\sum_{\alpha}\sigma(c_{\alpha l}\ot c_{jk})c_{i\alpha}\\
&=&\sum_{\alpha}x_{kl}^{j\alpha}c_{i\alpha}
\end{eqnarray*}
Hence
$$
\sum \sigma(c_{(1)}\ot h_{(1)})h_{(2)}c_{(2)}-
\sum \sigma(c_{(2)}\ot h)hc_{(1)}=\chi(i,j,k,l)=0
$$
i.e. (H1) also holds.

2. Suppose now that $R$ is bijective and let $S=R^{-1}$. Let $(y_{uv}^{ji})$
be a family of scalars of $k$ such that
$$
S(m_v\ot m_u)=\sum_{i,j}y_{uv}^{ji}m_i\ot m_j,
$$
for all $u$, $v=1,\cdots ,n$. As $RS=SR=Id_{M\ot M}$ we have
$$
\sum_{\alpha, \beta}x_{\beta \alpha}^{ip}y_{jq}^{\beta\alpha}=
\delta_{ij}\delta_{pq}, \qquad
\sum_{\alpha, \beta}y_{\beta \alpha}^{ip}x_{jq}^{\beta\alpha}=
\delta_{ij}\delta_{pq}
$$
for all $i$, $j$, $p$, $q=1,\cdots, n$. We define
$$
\sigma_{0}^{\prime}:C\ot C\to k, \quad
\sigma_{0}^{\prime}(c_{iv}\ot c_{ju}):=y_{uv}^{ji}
$$
for all $i$, $j$, $u$, $v=1,\cdots, n$. Now we extend $\sigma_{0}^{\prime}$
to a map $\sigma_{1}^{\prime}:C\ot T(C)\to k$ in such a way that
$\sigma_{1}^{\prime}$ satisfies (H2) and (H3). First we prove that
$\sigma_{1}^{\prime}$ is an inverse in the convolution algebra
$\Hom_k(C\ot T(C),k)$ of $\sigma_{1}$. Let $p$, $q$, $i$, $j=1,\cdots, n$.
We have:
\begin{eqnarray*}
\sum \sigma_1\Bigl((c_{pq})_{(1)}\ot (c_{ij})_{(1)}\Bigl)
\sigma_{1}^{\prime}\Bigl((c_{pq})_{(2)}\ot (c_{ij})_{(2)}\Bigl)&=&
\sum_{\alpha, \beta} \sigma_1(c_{p\alpha}\ot c_{i\beta})
\sigma_{1}^{\prime}(c_{\alpha q}\ot c_{\beta j})\\
&=&\sum_{\alpha, \beta}x_{\beta \alpha}^{ip}y_{jq}^{\beta\alpha}=
\delta_{ij}\delta_{pq}=\varepsilon(c_{ij})\varepsilon(c_{pq})
\end{eqnarray*}
and
\begin{eqnarray*}
\sum \sigma_1^{\prime}\Bigl((c_{pq})_{(1)}\ot (c_{ij})_{(1)}\Bigl)
\sigma_{1}\Bigl((c_{pq})_{(2)}\ot (c_{ij})_{(2)}\Bigl)&=&
\sum_{\alpha, \beta} \sigma_1^{\prime}(c_{p\alpha}\ot c_{i\beta})
\sigma_{1}(c_{\alpha q}\ot c_{\beta j})\\
&=&\sum_{\alpha, \beta}y_{\beta \alpha}^{ip}x_{jq}^{\beta\alpha}=
\delta_{ij}\delta_{pq}=\varepsilon(c_{ij})\varepsilon(c_{pq}).
\end{eqnarray*}
Hence, $\sigma_1\in \Hom_k(C\ot T(C), k)$ is invertible in convolution.
In order to prove that $\sigma\in \Hom_k(C\ot B(R), k)$ remains invertible
in the convolution, it is enought to prove that $\sigma_{1}^{\prime}$
factorizes to a map $\sigma^{\prime}:C\ot B(R)\to k$. We will prove now the
following:

{\sl $\sigma_{1}^{\prime}:C\ot T(C)\to k$ factorizes to a map
$\sigma^{\prime}:C\ot B(R)\to k$ if and only
if $R^{12}R^{13}=R^{13}R^{12}$.}

Indeed, $\sigma_{1}^{\prime}$ factorizes to a map
$\sigma^{\prime}:C\ot B(R)\to k$ if and only if, for any $i$, $j$, $k$,
$l$, $p$, $q=1,\cdots,n$, we have
$$
\sigma_{1}^{\prime}(c_{pq}\ot \chi(i,j,l,k))=0,
$$
which means
$$
\sum_{u,v}x_{uv}^{ji}\sigma_1^{\prime}(c_{pq}\ot c_{uk}c_{vl})=
\sum_{\alpha}x_{kl}^{j\alpha}\sigma_1^{\prime}(c_{pq}\ot c_{i\alpha})
$$
which is equivalent to
$$
\sum_{u,v,\beta}x_{uv}^{ji}\sigma_1^{\prime}(c_{p\beta}\ot c_{uk})
\sigma_1^{\prime}(c_{\beta q}\ot c_{vl})=
\sum_{\alpha}x_{kl}^{j\alpha}y_{\alpha q}^{ip}
$$
i.e.
$$
\sum_{u,v,\beta}x_{uv}^{ji}y_{k\beta}^{up}y_{lq}^{v\beta}=
\sum_{\alpha}x_{kl}^{j\alpha}y_{\alpha q}^{ip}.
$$
Now, from lemma (\ref{hei}) this equation is equivalent to
$$
R^{23}S^{13}S^{12}=S^{12}R^{23}.
$$
But $S=R^{-1}$, so the last equation turns into
\begin{equation}\label{adi}
R^{12}R^{23}=R^{23}R^{12}R^{13}.
\end{equation}
$R$ is a bijective solution of the Hopf equation, i.e.
$R^{12}R^{23}=R^{23}R^{13}R^{12}$. Hence, the equation (\ref{adi})
holds if and only if $R^{12}R^{13}=R^{13}R^{12}$.

This completes the proof of the theorem.
\end{proof}

\begin{remark}
There exists a major difference between the second point of our theorem and
the corresponding case for the quantum Yang-Baxter equation.
In this latter case, if $R$ is a bijective solution for the quantum
Yang-Baxter equation, then the map $\sigma :A(R)\ot A(R)\to k$,
which makes the bialgebra $A(R)$ co-quasitriangular, is invertible in
convolution.
Behind this lies the elementary observation that, if $R$ is a solution for
the quantum Young-Baxter equation, then $R^{-1}$ is also a solution. Now,
if $R$ is a bijective solution of the Hopf equation, then $R^{-1}=S$ is
not a solution for the Hopf equation; more precisely, $S$ is a solution for
the equation $S^{12}S^{13}S^{23}=S^{23}S^{12}$.
\end{remark}

Now we shall apply our theorem in order to give more examples of Hopf
functions on bialgebras.

\begin{examples}\label{210}
1. Let $q$ be a scalar of $k$ and $f_q:k^2\to k^2$, $f_q((x,y)):=(x+qy,0)$
for all $(x,y)\in k^2$, i.e., if $k=\nrea$, $f_q$ sends all the points of
the euclidian plane on the $Ox$ coordinate axis under an angle
arctg($q$) with respect to the $Oy$ axis. In \cite{M1} we constructed three
bialgebras associated to this map:
\begin{enumerate}
\item the first one, $B_{q}^{2}(k)$, corresponds to the solution of the Hopf
equation $f_q\ot (Id_{k^2}-f_q)$. $B_{q}^{2}(k)$ can be described as
follows:
\begin{enumerate}
\item as an algebra, $B_{q}^{2}(k)$ is generated by $x$, $y$, $z$ with
the relations
$$
yx=x,\qquad yz+qy^2=qx
$$
\item the comultiplication $\Delta$ and the counity $\varepsilon$
are given by:
\begin{equation}\label{ab}
\Delta(x)=x\ot x,\quad \Delta(y)=y\ot y,\quad \Delta(z)=x\ot z+z\ot y
\end{equation}
and
\begin{equation}\label{ac}
\varepsilon(x)=\varepsilon(y)=1,\quad \varepsilon(z)=0.
\end{equation}
\end{enumerate}
\item the second, $D_{q}^{2}(k)$, corresponds to $f_q\ot Id_{k^2}$ and
can be described as follows:
\begin{enumerate}
\item as an algebra, $D_{q}^{2}(k)$ is generated by $x$, $y$, $z$ with
the relations
$$
x^2=x=yx, \qquad zx=0, \qquad z^2+qzy=0, \qquad
xz+qxy=yz+qy^2=qx.
$$
\item the comultiplication $\Delta$ and the counity $\varepsilon$ are
given by equations (\ref{ab}) and (\ref{ac}).
\end{enumerate}
\item the third, $E_{q}^{2}(k)$, corresponds to $f_q\ot f_q$ and
can be described as follows:
\begin{enumerate}
\item as an algebra, $E_{q}^{2}(k)$ is generated by $x$, $y$, $z$ with
the relations
$$
x^2=x,\qquad xz+qxy=qx,\qquad zx+qyx=qx,\qquad
z^2+qyz+qzy+q^2y^2=q^2x.
$$
\item the comultiplication $\Delta$ and the counity $\varepsilon$ are
given by equations (\ref{ab}) and (\ref{ac}).
\end{enumerate}
\end{enumerate}
Let $C$ be the three dimensional subcoalgebra of $B_{q}^{2}(k)$
(resp. $D_{q}^{2}(k)$, $E_{q}^{2}(k)$) with $\{x,y,z\}$ a $k$-basis.
Then:
\begin{enumerate}
\item there exists a Hopf function
$$
\sigma :C\ot B_{q}^{2}(k)\to k
$$
such that
$$
\begin{array}{llll}
\sigma(x\ot 1)=1, &\sigma(x\ot x)=0, & \sigma(x\ot y)=1, & \sigma(x\ot z)=-q, \\
\sigma(y\ot 1)=1, &\sigma(y\ot x)=0, & \sigma(y\ot y)=0, & \sigma(y\ot z)=0,  \\
\sigma(z\ot 1)=0, &\sigma(z\ot x)=0, & \sigma(z\ot y)=q, & \sigma(z\ot z)=-q^2.
\end{array}
$$
\item there exists a Hopf function
$$
\sigma :C\ot D_{q}^{2}(k)\to k
$$
such that
$$
\begin{array}{llll}
\sigma(x\ot 1)=1, &\sigma(x\ot x)=1, & \sigma(x\ot y)=1, & \sigma(x\ot z)=0,\\
\sigma(y\ot 1)=1, &\sigma(y\ot x)=0, & \sigma(y\ot y)=0, & \sigma(y\ot z)=0,\\
\sigma(z\ot 1)=0, &\sigma(z\ot x)=q, & \sigma(z\ot y)=q, & \sigma(z\ot z)=0.
\end{array}
$$
\item there exists a Hopf function
$$
\sigma :C\ot E_{q}^{2}(k)\to k
$$
such that
$$
\begin{array}{llll}
\sigma(x\ot 1)=1, &\sigma(x\ot x)=1, & \sigma(x\ot y)=0, & \sigma(x\ot z)=q,\\
\sigma(y\ot 1)=1, &\sigma(y\ot x)=0, & \sigma(y\ot y)=0, & \sigma(y\ot z)=0,\\
\sigma(z\ot 1)=0, &\sigma(z\ot x)=q, & \sigma(z\ot y)=0, & \sigma(z\ot z)=q^2.
\end{array}
$$
\end{enumerate}
2. Let $k$ be a field of characteristic two and ${\cal F}(k)$ the five
dimensional noncommutative and noncocomutative bialgebra constructed in
\cite{M1},  corresponding to
$$
R=
\left(
\begin{array}{cccc}
1&0&0&0\\
0&1&1&0\\
0&0&1&0\\
0&0&0&1
\end{array}
\right)
$$
which is a solution of the Hopf equation. ${\cal F}(k)$ is described by the
following:

$\bullet$ as a vector space, ${\cal F}(k)$ is five dimensional with
$\{1,x,y,z,t \}$ a $k$-basis.

$\bullet$ the multiplication rule is given by:
$$
x^2=x, \quad y^2=z^2=0, \quad t^2=t,
$$
$$
xy=y,\quad yx=0,\quad xz=zx=z,\quad xt=t,\quad tx=x,
$$
$$
yz=0,\quad zy=x+t,\quad yt=0,\quad ty=y, \quad zt=tz=z.
$$
$\bullet$ the comultiplcation $\Delta$ and the counity $\varepsilon$
are given in such way that the matrix
$$
\left(
\begin{array}{cc}
x&y\\
z&t
\end{array}
\right)
$$
is comultiplicative. Let $C$ be the four dimensional subcoalgebra of
${\cal F}(k)$ with $\{x,y,z,t \}$ a $k$-basis. Then there exists a Hopf
function
$$
\sigma :C\ot {\cal F}(k)\to k
$$
such that
$$
\begin{array}{lllll}
\sigma(x\ot 1)=1, &\sigma(x\ot x)=1, & \sigma(x\ot y)=0,
& \sigma(x\ot z)=0, &\sigma(x\ot t)=1, \\
\sigma(y\ot 1)=0, &\sigma(y\ot x)=0, & \sigma(y\ot y)=0,
& \sigma(y\ot z)=1, &\sigma(y\ot t)=0, \\
\sigma(z\ot 1)=0, &\sigma(z\ot x)=0, & \sigma(z\ot y)=0,
& \sigma(z\ot z)=1, &\sigma(z\ot t)=0, \\
\sigma(t\ot 1)=1, &\sigma(t\ot x)=1, & \sigma(t\ot y)=0,
& \sigma(t\ot z)=0, &\sigma(t\ot t)=0,
\end{array}
$$
As $R$ is bijective and
$$
R^{12}R^{13}=R^{13}R^{12}=
\left(
\begin{array}{cccccccc}
1&0&0&0&0&0&0&0\\
0&1&0&0&1&0&0&0\\
0&0&1&0&1&0&0&0\\
0&0&0&1&0&1&1&0\\
0&0&0&0&1&0&0&0\\
0&0&0&0&0&1&0&0\\
0&0&0&0&0&0&1&0\\
0&0&0&0&0&0&0&1
\end{array}
\right)
$$
we obtain that $\sigma $ is invertible in convolution. Since $k$ has the
characteristic two, $R^{-1}=R$, hence $\sigma^{-1}=\sigma$.
\end{examples}

\section{Appendix: the Hopf elements}
In the preceding section we introduced the concept of bialgebra
$(H,C,\sigma)$ with a Hopf function $\sigma :C\ot H\to k$,
which, for the Hopf equation, plays the same role as the
co-quasitriangular bialgebra for the quantum Yang-Baxter equation. Now we
shall define the correspondent of the concept of quasitriangular bialgebra,
which is also involved in the quantum Yang-Baxter equation (see \cite{Dr}).
For the reasons presented in the preceding section, the definition of this
concept relative to a subalgebra becomes mandatory.
\begin{definition}
Let $H$ be a bialgebra and $A$ be a subalgebra of $H$. An element
$R=\sum R^1\ot R^2 \in A\ot H$ is called a Hopf element if:

$(HE \:1)\quad \sum \Delta(R^1)\ot R^2=R^{13}R^{23}$

$(HE \:2)\quad \sum \varepsilon(R^1)R^2=1$

$(HE \:3) \quad \Delta^{cop}(a)R=R(1\ot a)$

for all $a\in A$. In this case, we shall say that $(H,A,R)$ is a bialgebra
with  a Hopf element.
\end{definition}
\begin{remarks}
1. (HE 1) and (HE 2) are (QT 1) and (QT 2), respecting the definition
relative to the subalgebra $A$. (HE 3) is obtained by modifying the right
hand side of (QT 5), such that an integral type condition is obtained.
We shall detail:

Let $t:=\sum R^1\varepsilon(R^2)\in A$. (HE 3) can be written:
\begin{equation}\label{100}
\sum a_{(2)}R^1\ot a_{(1)}R^2=\sum R^1\ot R^2a
\end{equation}
for all $a\in A$. Applying $I\ot \varepsilon$ in this equation, we get
$$
at=\varepsilon(a)t
$$
for all $a\in A$. Hence, if $A$ is a subbialgebra of $H$, then $t$ is a left
integral in $A$. Now, if we apply $I\ot \varepsilon\ot \varepsilon$
to (HE 1) we get $t^2=t$. It follows that $t=tt=\varepsilon(t)t$, hence
$\varepsilon(t)=1$. Using the Maschke theorem for Hopf algebras, we
conclude that: if $(H,A,R)$ is a bialgebra with a Hopf element and
$A$ is a finite dimensional subbialgebra of $H$ with an antipode,
then $A$ is semisimple.

Conversely, if $t$ is a left integral in $A$, then $R:=t\ot 1$ satisfies
(HE 3).

2. Let $H$ be a bialgebra and $A$ be a subalgebra of $H$. Then $R=1\ot 1$
is a Hopf element if and only if $A=k$.

Indeed, if $R=1\ot 1$ then (HE 3) becomes $\Delta^{cop}(a)=1\ot a$, for all
$a\in A$. Hence, $a=\varepsilon(a)1_H$, for all $a\in A$, i.e. $A=k$.

3. Let $(H,A,R)$ be a bialgebra with a Hopf element. Suppose that $H$ has an
antipode $S$. Then $R$ is invertible and $R^{-1}=\sum S(R^1)\ot R^2$.
Moreover, if we denote $u:=\sum S(R^2)R^1\in H$ then
$$
S(a)u=\varepsilon(a)u,
$$
for all $a\in A$. This formula is obtained if we apply
$m_H\tau(I\ot S)$ in the equation (\ref{100}). We observe that if
$A\neq k$ then $u$ is not invertible (if $u$ is invertible, then
$R^{-1}=\sum S(R^1)\ot R^2=\sum \varepsilon(R^1)\ot R^2=1\ot 1$, i.e. $A=k$).
\end{remarks}

\begin{proposition}
Let $(H,A,R)$ be a bialgebra with a Hopf element $R\in A\ot H$. Then:
\begin{enumerate}
\item in the tensor product algebra $A\ot H\ot H$, the following identity
holds
\begin{equation}\label{101}
R^{23}R^{13}R^{12}=R^{12}R^{23}
\end{equation}
\item if $(M,\cdot)$ is a left $H$-module, then the map
$$
{\cal R}:M\ot M\to M\ot M, \quad
{\cal R}(m\ot n)=\sum R^1\cdot m\ot R^2\cdot n
$$
is a solution of the Hopf equation.
\end{enumerate}
\end{proposition}

\begin{proof}
1. (HE 1) is equivalent to

$(HE \:1^{\prime})\quad \sum \Delta^{cop}(R^1)\ot R^2=R^{23}R^{13}$

Now, for $r=R$, we have:
\begin{eqnarray*}
R^{23}R^{13}R^{12}&=&\sum \Bigl(\Delta^{cop}(R^1)\ot R^2 \Bigl)
(r^1\ot r^2\ot 1)\\
&=&\sum \Delta^{cop}(R^1)R\ot R^2\\
&=&\sum R(1\ot R^1)\ot R^2\\
&=&\sum r^1\ot r^2R^1\ot R^2\\
&=&R^{12}R^{23}
\end{eqnarray*}
2. Follows from equation (\ref{101}), as
$$
{\cal R}^{23}{\cal R}^{13}{\cal R}^{12}(l\ot m\ot n)=
R^{23}R^{13}R^{12}\cdot (l\ot m\ot n)
$$
and
$$
{\cal R}^{12}{\cal R}^{23}(l\ot m\ot n)=
R^{12}R^{23}\cdot (l\ot m\ot n)
$$
for all $l$, $m$, $n\in M$.
\end{proof}

Before presenting a few examples, we note that, if the compatibility
condition (HE 3) holds for a set of generators of $A$ as an algebra,
then it holds for any $a\in A$.

\begin{examples}
1. Let $H={\cal T}(k)$ be the three dimensional bialgebra from the
previous section and $A$ be the two dimensional subalgebra of ${\cal T}(k)$
generated by $x$. Let $R:=x\ot 1$. Then $R$ is a Hopf element.

Indeed, it is enough to check that (HE 3) holds for $a=x$. We have:
$$
\Delta^{cop}(x)R=(x\ot x)(x\ot 1)=x^2\ot x=x\ot x=(x\ot 1)(1\ot x),
$$
hence, $R$ is a Hopf element.

2. Let $H=E_q^2(k)$ and $A$ be the two dimensional subalgebra of
$E_q^2(k)$ generated by $x$.
Then $R:=x\ot 1$ is a Hopf element. The proof is the same as in the
previous example.

3. If in the above examples we have constructed bialgebras with Hopf elements
in which $A$ is finite dimensional, now we shall give an example in which
$A$ is infinite dimensional.

Let $H=D_q^2(k)$ and $A$ be the subalgebra of $D_q^2(k)$ generated by
$x$ and $y$. Then $R:=x\ot 1$ is a Hopf element.
We only verify for $a=y$ that (HE 3) holds (for $a=x$, the proof is the one
presented in the previous example). We have:
$$
\Delta^{cop}(y)(x\ot 1)=(y\ot y)(x\ot 1)=yx\ot y=x\ot y=(x\ot 1)(1\ot y)
$$
i.e. (HE 3) holds for $y$.
\end{examples}
%%%%%%%%%%%%%%%%%%%%%%%%%%%%%%%%%%%%%%%%%


\begin{thebibliography}{10}
\bibitem{CMZ1}
S. Caenepeel, G. Militaru, S. Zhu, Crossed modules and Doi-Hopf modules,
{\sl Israel J. Math.}, in press.
\bibitem{CMZ2}
S. Caenepeel, G. Militaru, S. Zhu, Doi-Hopf modules, Yetter-Drinfel'd modules
and Frobenius type properties, {\sl Trans. AMS}, in press.
\bibitem{CMZ3}
S. Caenepeel, G. Militaru, S. Zhu, A Maschke type theorem for Doi-Hopf
modules, {\sl J. Algebra} {\bf 187} (1997), 388-412.
\bibitem{D}
Y. Doi, Unifying Hopf modules, {\sl J. Algebra} {\bf 153} (1992), 373-385.
\bibitem{Dr}
V.G. Drinfel'd, Quantum groups, {\sl Proc. of the International Congress of
Mathematics}, 798-820, 1987.
\bibitem{FMS}
D. Fischman, S. Montgomery, H.-J. Schneider, Frobenius extensions of
subalgebras of Hopf algebras, {\sl Trans. AMS}, in press.
\bibitem{K}
C. Kassel, Quantum groups, Springer Verlag, Berlin, 1995.
\bibitem{LR}
L.A. Lambe, D. Radford, Algebraic aspects of the
quantum Yang-Baxter equation, {\sl J. Algebra} {\bf 54} (1992), 228-288.
\bibitem{MS}
G. Militaru, D. Stefan, Extending modules for Hopf Galois extensions,
{\sl Comm. Algebra} {\bf 14} (1994), 5657-5678.
\bibitem{M1}
G. Militaru, The Hopf modules category and nonlinear equations,
preprint 1997.
\bibitem{M}
S. Montgomery, Hopf algebras and their actions on rings, American Mathematical
Society, Providence, 1993.
\bibitem{R1}
D. Radford, Solutions to the quantum Yang-Baxter equation and the
Drinfel'd double, {\sl J. Algebra} {\bf 161} (1993), 20-32.
\bibitem{RT}
D. Radford, J. Towber, Yetter-Drinfel'd categories
associated to an arbitrary bialgebra, {\sl J. Pure and Appl. Algebra}
{\bf 87} (1993), 259-279.
\bibitem{Y}
D.N. Yetter, Quantum groups and representations of monoidal categories,
{\sl Math. Proc. Cambridge Philos. Soc.}, {\bf 108} (1990), 261-290.
\end{thebibliography}
\end{document}